\documentclass[11pt, twoside]{article}
\usepackage{amssymb, amsmath}
\usepackage{amscd}

\begin{document}
\newtheorem{Mthm}{Main Theorem.}
\newtheorem{Thm}{Theorem}[section]
\newtheorem{Prop}[Thm]{Proposition}
\newtheorem{Lem}[Thm]{Lemma}
\newtheorem{Cor}[Thm]{Corollary}
\newtheorem{Def}[Thm]{Definition}
\newtheorem{Guess}[Thm]{Conjecture}
\newtheorem{Ex}[Thm]{Example}
\newtheorem{Rmk}{Remark.}
\newtheorem{Not}{Notation.}
\newtheorem{thmA}{Theorem A.}
\newtheorem{thmB}{Theorem B.}
\newtheorem{thmC}{Theorem C.}
\newtheorem{DefProp}{Proposition-Definition}

\renewcommand{\theThm} {\thesection.\arabic{Thm}}
\renewcommand{\theProp}{\thesection.\arabic{Prop}}
\renewcommand{\theLem}{\thesection.\arabic{Lem}}
\renewcommand{\theCor}{\thesection.\arabic{Cor}}
\renewcommand{\theDef}{\thesection.\arabic{Def}}

\renewcommand{\theEx}{\thesection.\arabic{Ex}}
\renewcommand{\theDefProp}{\thesection.\arabic{DefProp}}
\renewcommand{\theRmk}{}
\renewcommand{\theMthm}{}
\renewcommand{\theNot}{}
\renewcommand{\thethmA}{}
\renewcommand{\thethmB}{}
\renewcommand{\thethmC}{}
\renewcommand{\thefootnote}{\fnsymbol{footnote}}

\newcommand{\spec}{\textnormal{Spec}\hspace{1mm}}
\newcommand{\pf}{{\bfseries\itshape Proof. }}
\newcommand{\mult}{\textnormal{mult}}
\newcommand{\com}{\hspace{-2mm}\textnormal{\textbf{.}}\hspace{2mm}}
\newcommand{\bir}{-\hspace{-1mm}\rightarrow}
\newcommand{\qed}{\hfill Q.E.D.\newline}
\newcommand{\lct}{\operatorname{lct}}
\newcommand{\tlct}{\operatorname{tlct}}

\begin{center}
\Large{\textbf{A NOTE ON DEL PEZZO FIBRATIONS OF DEGREE 1}}
\end{center}\vspace{15mm}
\begin{center}
\textbf{Jihun Park}
\vspace{5mm}\\
Department of Mathematics, The University of Georgia\\
Athens, GA 30602\\
E-mail : jhpark@math.uga.edu

\phantom{this is space} \vspace{10mm} \textbf{ABSTRACT}
\end{center}
The purpose of this paper is to extend the results in
\cite{P99} in the case of del Pezzo fibrations of degree 1. To
this end we investigate the anticanonical linear systems of del
Pezzo surfaces of degree 1. We then classify all possible
effective anticanonical divisors on Gorenstein del Pezzo surfaces
of degree 1 with canonical singularities. \vspace{3mm}

\thispagestyle{empty}
\renewcommand{\thesection}{\large{\arabic{section}.}}
\section{\hspace{-3mm}\large{Introduction}}
\renewcommand{\thesection}{\arabic{section}}
In \cite{P99}, it has been proven that smooth del Pezzo fibrations
of degree at most 4 over a discrete valuation ring cannot be
birationally transformed. Taking $\mathbb{P}^1$-bundle over a
discrete valuation ring into consideration, we find the result
quite interesting. Because the study of birational maps of del
Pezzo fibrations is quite important in birational geometry, we
need to pay attention to the result in \cite{P99}.

In the present paper, we extend the result in \cite{P99}.  The
result in \cite{P99} requires the special fibers to be smooth. It
is natural to expect the similar result when mild singularities
are allowed on the special fibers.

Let $\mathcal{O}$ be a discrete valuation ring with residue field
$k$ of characteristic 0. We denote by $K$ the quotient field of
$\mathcal{O}$. A model of a variety $X_K$ defined over $K$ is a
flat scheme $X$ defined over $\spec \mathcal{O}$ whose generic
fiber is isomorphic to $X_K$. Fano fibrations are models of smooth
Fano varieties defined over $K$. In particular, del Pezzo
fibrations of degree $d$ are models of smooth del Pezzo surfaces
of degree $d$ defined over $K$. We let $T=\spec \mathcal{O}$. For
a scheme $\pi : Z\longrightarrow T$, the scheme-theoretic fiber
$\pi^*(o)$ is denoted by $S_Z$, where $o$ is the closed point of
$T$.

Let $X/T$ be a Gorenstein del Pezzo
fibration. In the present paper, we always assume the following;
\begin{itemize}
\item
The special fiber $S_X$ is reduced and irreducible.
\item
Any birational map of $X$ into another del Pezzo fibration over
$T$ is biregular on generic fiber.
\end{itemize}
 Such del Pezzo fibrations are studied in \cite{Co96} and
\cite{Ko97a}. It should be remarked that the second condition
automatically holds if the generic fiber is a del Pezzo surface of
degree 1 with Picard number 1 (\cite{I96}). The author
(\cite{P99}) studied birational maps between them with one more
condition:
\begin{itemize}
\item
The log pair $(X, S_X)$ is purely log terminal.
\end{itemize}
We will also assume this condition throughout the paper.
The first and the last conditions force the special fiber
$S_X$ to be a normal Gorenstein del Pezzo surface
with canonical singularities.

\begin{Rmk}
\textnormal{In \cite{P99} we assume that del Pezzo fibrations are
$\mathbb{Q}$-factorial. But it turns out that we do not have to assume
$\mathbb{Q}$-factoriality. Indeed, it is enough to assume
only Gorenstein singularities.}
\end{Rmk}

In this paper, we are mainly interested in del Pezzo fibrations of
degree 1. The main task is to extend the results in \cite{P99}.
\begin{thmA}
Let $X/T$ and $Y/T$ be del Pezzo fibrations of degree $1$. We
suppose that the special fiber $S_X$ of $X$ has only singularities
of type $\textnormal{A}_n$, $n\geq 3$. Suppose that there is a
birational map $f:X\bir Y$. Then exactly one of the following holds:
\begin{itemize}
\item
The birational map $f$ is biregular.
\item
The special fiber $S_Y$ of $Y$
has a singularity of type $\textnormal{E}_8$. Moreover, it is a unique
singular point of $S_Y$.
\end{itemize}
Furthermore, if the anticanonical linear system $|-K_{S_X}|$ of
$S_X$ has no cuspidal rational curve, then only the first can
occur.
\end{thmA}
\begin{thmB}
Notation is as above. Suppose that the special fiber $S_X$ has
only singularities of type $\textnormal{A}_n$, $n\geq 1$. We
assume that $S_X$ has a singular point $p$ of type either
$\textnormal{A}_1$ or $\textnormal{A}_2$ and  that $f$ is not
biregular.
\begin{enumerate}
\item
There is a cuspidal rational curve $D$ in $|-K_{S_X}|$.
\item
If every cuspidal rational curve in $|-K_{S_X}|$ has a cusp
outside of singular points of type $\textnormal{A}_2$ and $\textnormal{A}_1$,
then $S_Y$ has only one singular point and the
singularity type is $\textnormal{E}_8$.
\end{enumerate}
Suppose that there is an element in $|-K_{S_X}|$ which has a cusp
at the point $p$.
\begin{enumerate}
\item[3.]
If the point $p$ has singularity type $\textnormal{A}_1$ and every
element of $|-K_{S_X}|$ has no cusp at  singular points different
type from $\textnormal{A}_1$,
then $S_Y$ has a singular point of type either
$\textnormal{E}_8$ or $\textnormal{E}_7$.
\item[4.]
If the point $p$ has singularity type $\textnormal{A}_2$, then
$S_Y$ has a singular point of type $\textnormal{E}_8$,
$\textnormal{E}_7$, or $\textnormal{E}_6$.
\end{enumerate}
\end{thmB}
\begin{thmC}
Notation is as above. Suppose that the special fiber $S_X$ has
only singularities of type $\textnormal{A}_n, n\geq 1$.
\begin{enumerate}
\item
If $S_Y$ has a singularity of type $\textnormal{E}_6$, then
$S_X$ has a cuspidal rational curve in $|-K_{S_X}|$ the cusp of
which is an $\textnormal{A}_2$-singularity of $S_X$.
\item
If $S_Y$ has a singularity of type $\textnormal{E}_7$, then $S_X$
has a cuspidal rational curve in $|-K_{S_X}|$ the cusp of which
is a singular point of $S_X$ of type either $\textnormal{A}_1$ or
$\textnormal{A}_2$.
\end{enumerate}
\end{thmC}

Before we proceed, we provide an example which
illustrates Theorem~A.
Let $X$ and $Y$ be  subschemes of
$\mathbb{P}^3_{\mathcal{O}}(1,1,2,3)$ defined by equations
$w^2+z^3+xy^5+x^5y=0$ and $w^2+z^3+xy^5+t^{24}x^5y=0$,
respectively, where $z$ and $w$ are of weight 2 and 3,
respectively, and $t$ is a local parameter of $\mathcal{O}$.
The special fiber of $X$ is smooth. On the other hand,
that of $Y$ has a single singular point which is of type $\textnormal{E}_8$.
There
is a birational map $f:X\bir Y$ defined by
$f(x,y,z,w)=(x,t^6y,t^{10}z,t^{15}w)$.

\renewcommand{\thesection}{\large{\arabic{section}.}}
\section{\hspace{-3mm}\large{Definitions}}
\renewcommand{\thesection}{\arabic{section}}
In the present section, we briefly explain the essential
definitions for this paper. Basically,  we use the definitions in
birational geometry ``textbooks". For the definitions of log
canonical singularities, pure log canonical singularities, centers
of log canonicity, 1-complements and so forth (see \cite{Koetal}
and \cite{Sho92}).

The most important tool in this paper is the concept of log
canonical threshold which was introduced by V.~V.~Shokurov.
\begin{Def}\com
Let $(X, B)$ be a log canonical pair. Let $D$ be a
$\mathbb{Q}$-Cartier divisor on $X$. The log canonical threshold of $D$
on $(X,B)$ is the number
\[\lct(D;X,B)=\operatorname{max}\{c : K_X+B+cD \textnormal{ is log
canonical}\}.\]
\end{Def}
It is easy to check that $0\leq \lct(D;X,B)\leq 1$.

We need to investigate anticanonical linear systems to obtain our
main results. To this end we will use the following number
to measure how bad an effective anticanonical divisor can be.
\begin{Def}\com
Let $X$ be a normal variety with nonempty anticanonical linear
system. We suppose that the log pair $(X,0)$ is log canonical. The
total log canonical threshold of $X$ is the number
\[\tlct(X)=\operatorname{max}\{c : K_X+cD \textnormal{ is log
canonical for any } D\in|-K_X|\}.\]
\end{Def}
Note that $0\leq \tlct(X)\leq 1$.

In the study of smooth weak del Pezzo surfaces, fundamental
cycles, introduced by M.~Artin (\cite{Ar66}),
give us a great deal of information. In \cite{Rei97}, they are
called numerical cycles.
\begin{DefProp}\com
Let $\pi:Y\longrightarrow X$ be a resolution of a point $p$ on a
normal surface $X$. Let $E=\sum E_i$ be the divisor of
the $\pi$-exceptional locus. Then there exists a unique effective
exceptional divisor $\Gamma=\sum a_iE_i$ such that $\Gamma>0$,
$\Gamma\cdot E_i\leq 0$ for every $E_i$, and $\Gamma$ is minimal with
respect to this property. The divisor $\Gamma$ is called the
fundamental cycle of
the bunch $\{E_i\}$.
\end{DefProp}
\pf See \cite{Rei97}. \qed

For a minimal resolution of a Du Val singularity, we can easily
find the corresponding fundamental cycle. Since the fundamental
cycles related to Du Val singularities are essential in this work
we list all of them in the Appendix. We
will see Kodaira's classification of degenerations of elliptic
curves in the Appendix.
\renewcommand{\thesection}{\large{\arabic{section}.}}
\section{\hspace{-3mm}\large{Normal Gorenstein del Pezzo surfaces
of degree 1 with canonical singularities}}
\renewcommand{\thesection}{\arabic{section}}
Let $S$ be a normal Gorenstein del Pezzo surface of degree 1 with
canonical singularities. And let $\pi:\tilde{S}\longrightarrow S$
be the minimal resolution of $S$. Then the smooth surface
$\tilde{S}$ is usually
called a weak del Pezzo surface because $-K_{\tilde{S}}$ is nef
and big. In the present section we study normal Gorenstein del Pezzo
surfaces of degree 1 with canonical singularities via smooth weak
del Pezzo surfaces. These surfaces were investigated
in \cite{Dema80}, \cite{Furu86},
\cite{HiWa81}, \cite{MiZha88}, \cite{MiZha93},
\cite{U83}, and
\cite{Zha88}.
In our study we pay more attention
to their effective anticanonical divisors.

\begin{Lem}\label{comp}\com
Any element of $|-K_S|$ is reduced and irreducible.
\end{Lem}
\pf Let $D=\sum a_i D_i$, where each $D_i$ is a prime divisor.
Since the surface $S$ is Gorenstein, each $a_i$ is a nonnegative integer. Then
\[1=D\cdot (-K_S)=\pi^*(D)\cdot \pi^*(-K_S)=\sum a_i\tilde{D}_i\cdot
(-K_{\tilde{S}})\geq \sum a_i,\] where $\tilde{D_i}$ is the
strict transform of $D_i$ via $\pi$. Note that no $\tilde{D_i}$ is
a $-2$-curve. Therefore, $D$ is reduced and irreducible. \qed
\begin{Lem}\label{Dema}\com
Let $H$ be an element of $|-K_{\tilde{S}}|$. And let $\Gamma$ be a
fundamental cycle of the minimal resolution
$\pi:\tilde{S}\longrightarrow S$. If $H$ contains a point of
$\Gamma$, then $H=\tilde{D}+\Gamma$, where $\tilde{D}$ is a
$-1$-curve.
\end{Lem}
\pf See \cite[p.53, Corollaire 2]{Dema80}.\qed

Let $D$ be an element of $|-K_S|$. We consider the pull-back
$\pi^*(D)$ of $D$ via $\pi$. Then we may write
\[\pi^*(D)=\tilde{D}+E,\]
where $\tilde{D}$ is the strict transform of $D$ and $E$ consists
of $-2$-curves.

\begin{Thm}\label{lc}\com
Let $p$ be a singular point of $S$. Suppose that $D$ passes
through the point $p$. Then the point $p$ is of type
$\textnormal{A}_n, n\geq 3$ \textnormal{(}$\textnormal{D}_n$,
$\textnormal{E}_6$, $\textnormal{E}_7$, $\textnormal{E}_8$,
resp.\textnormal{)}
if and only if the log canonical threshold of $D$ with respect to
$K_S$ is $1$ \textnormal{(}$\frac{1}{2}$, $\frac{1}{3}$, $\frac{1}{4}$,
$\frac{1}{6}$, resp.\textnormal{)}. Moreover, the divisor $D$ has a node at
$p$ if and only if $K_S+D$ is log canonical.
\end{Thm}
\pf  Since the pull-back $\pi^*(D)$ of $D$ belongs to
$|-K_{\tilde{S}}|$, Lemma~\ref{Dema} implies
\[\pi^*(K_S+cD)=K_{\tilde{S}}+c\tilde{D}+c\Gamma,\]
where $c$ is a constant and $\Gamma$ is the fundamental cycle on $\tilde{S}$
associated to the point $p$. Therefore, it is enough to
consider the log canonical threshold of $\tilde{D}+\Gamma$ with respect to
$K_{\tilde{S}}$.

{\bfseries \itshape Claim.}
$\tilde{D}+\operatorname{supp}(\Gamma)$ is a simple normal
crossing divisor.\\
Let $\Gamma=\sum_{i=1}^{n} a_iE_i$, where each $E_i$ is a
$-2$-curve. Observe the equations
\[0=\pi^*(D)\cdot E_j=\tilde{D}\cdot E_j +\sum_{i=1}^n a_iE_i\cdot
E_j, \hspace{5mm} j=1,\cdots, n.\] Then we obtain
\[M\left(\begin{array}{c}a_1\\a_2\\\vdots\\a_{n-1}\\a_n\end{array}\right)=
-\left(\begin{array}{c}\tilde{D}\cdot E_1\\\tilde{D}\cdot E_2\\
\vdots\\\tilde{D}\cdot E_{n-1}\\
\tilde{D}\cdot E_n\end{array}\right),\] where $M=(E_i\cdot E_j)$.
Looking at the fundamental cycles case by case (see the Appendix
at the end of the paper), we then see that the claim is true.
Moreover the configuration of the effective anticanonical divisor
$\pi^*(D)$ has a perfect match to one of Kodaira's singular
elliptic fibers, $\widetilde{\textnormal{E}}_8
\textnormal{(II}^*\textnormal{)}$, $\widetilde{\textnormal{E}}_7
\textnormal{(III}^*\textnormal{)}$, $\widetilde{\textnormal{E}}_6
\textnormal{(IV}^*\textnormal{)}$, $\widetilde{\textnormal{D}}_m
\textnormal{(I}^*_{m-4}\textnormal{)}$, and
$\widetilde{\textnormal{A}}_n
\textnormal{(I}_{n+1}\textnormal{)}$, $n\geq 3$.

Since $\tilde{D}$ is a smooth curve, the first statement
immediately follows from the claim.

The Appendix shows that the effective anticanonical divisor $\pi^*(D)$
is a wheel if and only if the point $p$ is a singularity of type
$\textnormal{A}_n$. This implies the second statement.\qed

\begin{Thm}\label{A1}\com
Notation as in Theorem~\ref{lc}.
\begin{enumerate}
\item
The point $p$ is of type $\textnormal{A}_1$
$(\textnormal{A}_2$ resp.\textnormal{)}
and $D$ has a cusp at
$p$ if and only if the log canonical threshold of $D$ is
$\frac{3}{4}$ \textnormal{(}$\frac{2}{3}$ resp.\textnormal{)}.
\item
If the point $p$ is of type either $\textnormal{A}_1$
or $\textnormal{A}_2$ and $K_S+D$ is log
canonical, then $D$ has a node at the point $p$.
\end{enumerate}
\end{Thm}
\pf The proof is similar to that of Theorem~\ref{lc}. We prove
the case corresponding to $\textnormal{A}_1$. Note that
$\tilde{D}$ can meet $\Gamma$ transversally or tangentially with
intersection number 2. If $\tilde{D}$ meets $\Gamma$
transversally, then $K_S+D$ is log canonical and $D$ has a node
at the point $p$. If $\tilde{D}$ meets $\Gamma$ tangentially,
then we can easily check that the log canonical threshold of
$\tilde{D}$ with respect to $K_{\tilde{S}}$ is $\frac{3}{4}$. And
$D$ has a cusp at the point $p$.

In the case of $\textnormal{A}_2$, it follows from the Appendix that
$\tilde{D}$ meets $\Gamma=E_1+E_2$ at two distinct points with
$E_1\cdot\tilde{D}=E_2\cdot\tilde{D}=1$ or at single point with
$E_1\cdot\tilde{D}=E_2\cdot\tilde{D}=1$. The first gives us the
log canonical threshold 1 and the latter provides the log
canonical threshold with $\frac{2}{3}$. \qed
\begin{Rmk}
\textnormal{In the proof of Theorem~\ref{A1}, we can also find
perfect matches to Kodaira's singular elliptic fibers. Namely, the
singularity  type $\textnormal{A}_1$ corresponds to
$*\widetilde{\textnormal{A}}_1\textnormal{(III}\textnormal{)}$ or
$\widetilde{\textnormal{A}}_1 \textnormal{(I}_{2}\textnormal{)}$
and the singularity type $\textnormal{A}_2$ corresponds to
$*\widetilde{\textnormal{A}}_2\textnormal{(IV}\textnormal{)}$ or
$\widetilde{\textnormal{A}}_2 \textnormal{(I}_{3}\textnormal{)}.$}
\end{Rmk}
Summarizing these results, we obtain perfect information on
total log canonical thresholds of normal Gorenstein del Pezzo
surfaces of degree 1 with canonical singularities. We can see a
beautiful numerical connection between log canonical thresholds and
Kodaira's classification of degenerations of elliptic curves.
\begin{Cor}\label{totallc}\com
Let $S$ be a normal Gorenstein del Pezzo surface of degree $1$
with canonical singularity as before.
\begin{itemize}
\item
$\widetilde{\textnormal{E}}_8
\textnormal{(II}^*\textnormal{)}.$\\
$\tlct(S)=\frac{1}{6}$ if and only if $S$ has an
$\textnormal{E}_8$ singularity.
\item
$\widetilde{\textnormal{E}}_7
\textnormal{(III}^*\textnormal{)}.$\\
$\tlct(S)=\frac{1}{4}$ if and only if $S$ has an
$\textnormal{E}_7$ singularity and no $\textnormal{E}_8$
singularity.
\item
$\widetilde{\textnormal{E}}_6
\textnormal{(IV}^*\textnormal{)}.$\\
$\tlct(S)=\frac{1}{3}$ if and only if $S$ has an
$\textnormal{E}_6$ but neither $\textnormal{E}_7$  nor
$\textnormal{E}_8$.
\item
$\widetilde{\textnormal{D}}_n
\textnormal{(I}^*_{n-4}\textnormal{)}$, $4\leq n\leq 8$.\\
$\tlct(S)=\frac{1}{2}$ if and only if $S$ has a
$\textnormal{D}_n$ singularity and no exceptional type
singularity.
\item
$*\widetilde{\textnormal{A}}_2 \textnormal{(IV}\textnormal{)}.$\\
$\tlct(S)=\frac{2}{3}$ if and only if $S$ has only
$\textnormal{A}_n$ type singularities and there is an element in
$|-K_S|$ which has a cusp at an $\textnormal{A}_2$ singularity of
$S$.
\item
$*\widetilde{\textnormal{A}}_1\textnormal{(III}\textnormal{)}.$\\
$\tlct(S)=\frac{3}{4}$ if and only if $S$ has only
$\textnormal{A}_n$ type singularities and $|-K_S|$ contains an
element having a cusp at an $\textnormal{A}_1$ singularity of $S$
but no element having a cusp at an $\textnormal{A}_2$ singularity
of $S$.
\item
$*\widetilde{\textnormal{A}}_0
\textnormal{(II}\textnormal{)}.$\\
$\tlct(S)=\frac{5}{6}$ if and only if $S$ has only
$\textnormal{A}_n$ type singularities and $|-K_S|$ contains an
element having a cusp but no element having a cusp at a singular
point of $S$.
\item
$\widetilde{\textnormal{A}}_n
\textnormal{(I}_{n+1}\textnormal{)}$, $n\leq 8$.
\\ $\tlct(S)=1$ if and only
if $S$ has only $\textnormal{A}_n$ type singularities and there
is no cuspidal rational curve in $|-K_S|$.
\end{itemize}
\end{Cor}
\pf Let $D$ be an effective anticanonical divisor of $S$.
By Lemma~\ref{comp}, $D$ is irreducible and reduced.
If $D$ passes through a singular point of $S$, then
the log canonical threshold of $D$ should obey the rules
in Theorems~\ref{lc} and~\ref{A1}.

Suppose that $D$ does not pass through any singular point of $S$.
We can easily check that the log canonical threshold of $D$ is
either 1 or $\frac{5}{6}$. Moreover,  $D$ has a cusp if and only
if the log canonical threshold  of $D$ is $\frac{5}{6}$. \qed
\begin{Cor}\com
Every effective anticanonical divisor of $S$ passing through a
singular point of type $\textnormal{D}_n$, $\textnormal{E}_6$,
$\textnormal{E}_7$, or $\textnormal{E}_8$ has a cusp at the singular point.
If it passes
through a singular point of type $\textnormal{A}_n$, $n\geq 1$,
then it has either a node or a cusp at the singular point.
\end{Cor}
\pf Let $D$ be an effective anticanonical divisor of $S$ passing
through a singular point. Note that $D$ is a Cartier divisor. It
then follows from Inversion of adjunction (\cite{Koetal}) that
$D$ cannot be smooth. The arithmetic genus of $D$ is 1 while that
of the strict transform of $D$ via $\pi$ is 0. The configuration
of $\pi^*(D)$ then implies the results. \qed

Slightly changing the focus, we look at the number of
singularities on the surface $S$ when it has an exceptional
singularity. This observation will later give some remarks on del Pezzo
fibrations. The following proposition can be easily derived
from \cite{MiZha88}.

\begin{Prop}\label{one}\com
Let $p$ be a singular point of $S$.
\begin{enumerate}
\item
The point $p$ is neither $\textnormal{A}_n$ nor $\textnormal{D}_n$,
$n\geq 9$.
\item
If $p$ is of type $\textnormal{A}_8$, $\textnormal{D}_8$, or
$\textnormal{E}_8$, then it is a
unique singular point on $S$.
\item
If $p$ is of type $\textnormal{A}_7$, $\textnormal{D}_7$, or
$\textnormal{E}_7$, then the surface $S$
has at most two singular points. Moreover, the possible
extra singularity is of type $\textnormal{A}_1$.
\item
If $p$ is of type $\textnormal{E}_6$, then the surface $S$ has at
most two singular points. The possible extra singularity is of
type either $\textnormal{A}_1$ or $\textnormal{A}_2$.
\end{enumerate}
\end{Prop}
\pf Since the minimal resolution $\tilde{S}$ of the surface $S$
has degree 1, the rank of $\tilde{S}$ is 9. Therefore, the number
of $-2$-curves is at most 8. This fact implies the first three
statements. For the last statement, we have to show that the
surface $S$ cannot have the three singularities $\textnormal{E}_6$,
$\textnormal{A}_1$ and $\textnormal{A}_1$ at the same time.
Because the rank of $S$
is 1, this can be verified by the classification of singularities
on normal Gorenstein del Pezzo surfaces of rank 1 in
\cite{MiZha88}. \qed

\renewcommand{\thesection}{\large{\arabic{section}.}}
\section{\hspace{-3mm}\large{Proofs of Theorems A, B, and C}}
\renewcommand{\thesection}{\arabic{section}}
For the readers' convenience we state the result in \cite{P99} which
is the main method for our proofs.
\begin{Thm}\label{P}\com
Notation is as in Theorem A. If the del Pezzo fibrations $X/T$
and $Y/T$ satisfy the conditions below, then $f$ is biregular.
\begin{itemize}
\item
\textnormal{(}Special fiber condition\textnormal{)}\\
The special fiber $S_X$ \textnormal{(}$S_Y$ resp.\textnormal{)}
is reduced and irreducible. And the log pair $(X, S_X)$
\textnormal{(}$S_Y$ resp.\textnormal{)} is purely log terminal.
\item
\textnormal{(}1-complement condition\textnormal{)}\\
For any $C\in|-K_{S_X}|$ \textnormal{(}$|-K_{S_Y}|$, resp.\textnormal{)},
there exists
$1$-complement $K_{S_X}+C_X$
\textnormal{(}$K_{S_Y}+S_Y$ resp.\textnormal{)} of $K_{S_X}$
\textnormal{(}$K_{S_Y}$, resp.\textnormal{)} such that $C_X$
\textnormal{(}$C_Y$, resp.\textnormal{)} does not contain
any center of log canonicity of $K_{S_X}+C$
\textnormal{(}$K_{S_Y}+C$ resp.\textnormal{)}.
\item
\textnormal{(}Surjectivity condition\textnormal{)}\\
Any $1$-complement of $K_{S_X}$ \textnormal{(}$K_{S_Y}$,
resp.\textnormal{)} can be extended to a $1$-complement of
$K_X+S_X$ \textnormal{(}$K_Y+S_Y$, resp.\textnormal{)}.
\item
\textnormal{(}Total lc threshold condition\textnormal{)}\\
The inequality $\tlct(S_X)+\tlct(S_Y)>1$ holds.
\end{itemize}
\end{Thm}
\pf See \cite{P99}\qed
In our situation, the first three conditions are satisfied (see
\cite{P99}). The main recipe for the proofs of Theorems A, B, and
C is the Total lc threshold condition.

{\bfseries \itshape Proof of Theorem A.} Suppose that the
birational map $f$ is not biregular. Then
$\tlct(S_X)+\tlct(S_Y)\leq 1$ by the Total lc threshold condition.
We also obtain $\tlct(S_X)\geq \frac{5}{6}$ from
Corollary~\ref{totallc}. Hence, applying Corollary~\ref{totallc}
and Proposition~\ref{one} to the special fiber $S_Y$ completes the
proof of the first assertion.

In order to prove the second assertion, it is enough to see that
$\tlct(S_X)=1$ and $\tlct(S_Y)\geq\frac{1}{6}$ under the
conditions. \qed

{\bfseries \itshape Proof of Theorem B.} For the first statement
suppose there is no cuspidal rational curve in $|-K_{S_X}|$. We
then obtain $\tlct(S_X)=1$. Since $\tlct(S_Y)$ is always positive,
this is a contradiction.

In the second case, Corollary~\ref{totallc} implies
$\tlct(S_X)=\frac{5}{6}$ because we assumed that $f$ is not
biregular. Therefore, $\tlct(S_Y)$ must be $\frac{1}{6}$. The
result then follows from Corollary~\ref{totallc}.

As for the third statement, we see that the conditions force
$\tlct(S_X)$ to be $\frac{3}{4}$. Therefore
$\tlct(S_Y)\leq\frac{1}{4}$. Then Corollary~\ref{totallc} shows
the result.

In the last statement,  the conditions give us
$\tlct(S_X)=\frac{2}{3}$. Thus the statement follows from
Corollary~\ref{totallc}.\qed
\begin{Rmk}
\textnormal{In the statements 3 and 4, we see more than the
existence of a certain type of singularity on the special fiber
$S_Y$. To be precise, Proposition~\ref{one} states the number and
types of singularities on $S_Y$ as in the case of
$\textnormal{E}_8$.}
\end{Rmk}

{\bfseries \itshape Proof of Theorem C.} Noting
Proposition~\ref{one}, we see that the conditions imply
$\tlct(S_Y)=\frac{1}{3}$ ($\frac{1}{4}$ resp.). Therefore
$\tlct(S_X)\leq\frac{2}{3}$ ($\frac{3}{4}$ resp.). Since $S_X$ has
only singularities of type $A_n$, the result follows from
Corollary~\ref{totallc}. \qed
\renewcommand{\thesection}{}
\section{\hspace{-3mm}\large{Appendix}}
\small Each diagram represents the fundamental cycle associated to
the given singularity type. These fundamental cycles
can be easily derived by simple computation. We may refer to
\cite{Rei97} or \cite{Ar66}. In each matrix equation the column
vector of the left hand side represents the multiplicities of the
$-2$-curves of the fundamental cycle in suitable order.
\begin{enumerate}
\item $\textnormal{A}_n$, $n\geq 1$
\begin{center}
\begin{picture}(500,40)(55,20)
\put(150,40){\circle*{5}}
\put(210,40){\circle*{5}} \put(180,40){\circle*{5}}
\put(150,40){\line(1,0){30}}
\put(210,40){\line(1,0){30}} \put(190,39){...}
\multiput(152,45)(30,0){4}{1}\put(240,40){\circle*{5}}
\end{picture}
\end{center}
\vspace{5mm}
\[\left(\begin{array}{cccccc}
        -2&1&0&\cdots&0&0 \\
        1&-2&1&0&\cdots&0 \\
        0&1&-2&1&0&\cdots \\
        \multicolumn{6}{c}\dotfill\\
        0&\cdots&0&1&-2&1 \\
        0&0&\cdots&0&1&-2
\end{array}\right)
\left(\begin{array}{c}1\\1\\1\\\vdots\\1\\1\end{array}\right)=
-\left(\begin{array}{c}1\\0\\0\\\vdots\\0\\
1\end{array}\right)
\]
where $n\geq 2$ for the matrix equation.
\item $\textnormal{D}_n$, $n\geq 4$
\begin{center}
\begin{picture}(500,40)(40,20)
\multiput(130,60)(100,0){1}{\circle*{5}}
\multiput(130,20)(100,0){1}{\circle*{5}} \put(150,40){\circle*{5}}
\put(210,40){\circle*{5}} \put(180,40){\circle*{5}}
\put(130,60){\line(1,-1){20}}
\put(130,20){\line(1,1){20}}
\put(150,40){\line(1,0){30}}
\put(210,40){\line(1,0){30}} \put(190,39){...}
\multiput(120,60)(82,0){1}{1} \multiput(120,18)(82,0){1}{1}
\multiput(149,45)(30,0){3}{2}\put(240,40){\circle*{5}}
\put(239,45){1}
\end{picture}
\end{center}
\vspace{5mm}
\[\left(\begin{array}{cccccccc}
        -2&1&0&0&0&\cdots&0&0 \\
        1&-2&1&0&0&\cdots&0&0 \\
        0&1&-2&1&0&\cdots&0&0 \\
        \multicolumn{8}{c}\dotfill\\
        0&\cdots&0&1&-2&1&0&0 \\
        0&0&\cdots&0&1&-2&1&1 \\
        0&0&0&\cdots&0&1&-2&0 \\
        0&0&0&\cdots&0&1&0&-2
\end{array}\right)
\left(\begin{array}{c}1\\2\\2\\\vdots\\2\\2\\1\\1\end{array}\right)=
-\left(\begin{array}{c}0\\1\\0\\0\\\vdots\\0\\0\\0
\end{array}\right)\]

\item $\textnormal{E}_6$
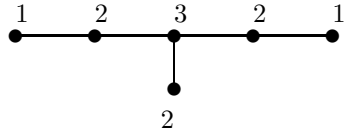
\begin{center}
\begin{picture}(500,20)(10,20)
\multiput(100,40)(30,0){5}{\circle*{5}} \put(160,20){\circle*{5}}
\put(100,40){\line(1,0){120}} \put(160,40){\line(0,-1){20}}
\put(100,45){1} \put(130,45){2} \put(160,45){3} \put(190,45){2}
\put(220,45){1}  \put(155,5){2}
\end{picture}
\end{center}
\vspace{5mm}
\[\left(\begin{array}{cccccc}
    -2&1&0&0&0&0 \\
    1&-2&1&0&0&0 \\
    0&1&-2&1&0&1 \\
    0&0&1&-2&1&0 \\
    0&0&0&1&-2&0 \\
    0&0&1&0&0&-2
    \end{array}\right)
\left(\begin{array}{c}1\\2\\3\\2\\1\\2\end{array}\right)=
-\left(\begin{array}{c}0\\0\\0\\0\\0\\
1\end{array}\right)
\]
\item $\textnormal{E}_7$
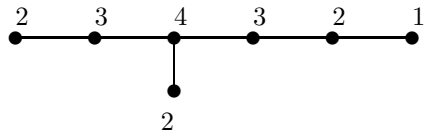
\begin{center}
\begin{picture}(500,20)(10,20)
\multiput(100,40)(30,0){6}{\circle*{5}} \put(160,20){\circle*{5}}
\put(100,40){\line(1,0){150}} \put(160,40){\line(0,-1){20}}
\put(100,45){2} \put(130,45){3} \put(160,45){4} \put(190,45){3}
\put(220,45){2} \put(250,45){1}  \put(155,5){2}
\end{picture}
\end{center}
\vspace{5mm}
\[\left(\begin{array}{ccccccc}
        -2&1&0&0&0&0&0 \\
        1&-2&1&0&0&0&0 \\
        0&1&-2&1&0&0&0 \\
        0&0&1&-2&1&0&1 \\
        0&0&0&1&-2&1&0 \\
        0&0&0&0&1&-2&0 \\
        0&0&0&1&0&0&-2
\end{array}\right)
\left(\begin{array}{c}1\\2\\3\\4\\3\\2\\2\end{array}\right)=
-\left(\begin{array}{c}0\\0\\0\\0\\
0\\1\\0\end{array}\right)\]
\item $\textnormal{E}_8$
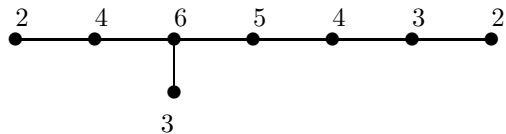
\begin{center}
\begin{picture}(500,20)(10,20)
\multiput(100,40)(30,0){7}{\circle*{5}} \put(160,20){\circle*{5}}
\put(100,40){\line(1,0){180}} \put(160,40){\line(0,-1){20}}
\put(100,45){2} \put(130,45){4} \put(160,45){6} \put(190,45){5}
\put(220,45){4} \put(250,45){3}  \put(155,5){3} \put(280,45){2}
\end{picture}
\end{center}
\vspace{5mm}
\[
\left(
\begin{array}{cccccccc}
    -2&1&0&0&0&0&0&0 \\
    1&-2&1&0&0&0&0&0 \\
    0&1&-2&1&0&0&0&0 \\
    0&0&1&-2&1&0&0&0 \\
    0&0&0&1&-2&1&0&1 \\
    0&0&0&0&1&-2&1&0 \\
    0&0&0&0&0&1&-2&0 \\
    0&0&0&0&1&0&0&-2
  \end{array}\right)
\left(\begin{array}{c}2\\3\\4\\5\\6\\4\\2\\3\end{array}\right)=
-\left(\begin{array}{c}1\\0\\0\\0\\0\\0\\0\\0
\end{array}\right)\]
\end{enumerate}

\textbf{Acknowledgement.} The author would like to express deep
gratitude to Professor Vyacheslav V. Shokurov, who gave him many
pieces of invaluable advice. This work was partially supported by
NSF Grant DMS-0100991.

\bibliographystyle{amsplain}

\providecommand{\bysame}{\leavevmode\hbox to3em{\hrulefill}\thinspace}

\end{document}